\documentclass[11pt]{amsart}

\usepackage{epsfig}
\usepackage{verbatim}
\usepackage{amssymb}
\usepackage{amsfonts}
\usepackage{amsbsy}

\newcommand{\ds}{\displaystyle}

\newcommand{\FT}{\widehat}

\newcommand{\eproof}{\hfill\rule{2.2mm}{3.0mm}}
\newcommand{\esubproof}{\hfill$\square$}
\newcommand{\Proof}{\noindent {\bf Proof.~~}}
\newcommand{\A}{{\mathcal A}}
\newcommand{\B}{{\mathcal B}}
\newcommand{\D}{{\mathcal D}}
\newcommand{\E}{{\mathcal E}}
\renewcommand{\SS}{{\mathcal S}}

\renewcommand{\P}{{\mathcal P}}
\newcommand{\QQ}{{\mathcal Q}}
\newcommand{\R}{{\mathbb R}}
\newcommand{\Z}{{\mathbb Z}}

\newcommand{\Q}{{\mathbb Q}}

\newcommand{\wmod}[1]{\mbox{~({\rm mod}~$#1$)}}
\renewcommand{\eqref}[1]{(\ref{#1})}
\newcommand{\shsp}{\hspace{1em}}
\newcommand{\mhsp}{\hspace{2em}}
\newcommand{\tmod}[1]{({\rm mod}~#1)}
\newtheorem{prop}{Proposition}[section]
\newtheorem{lem}[prop]{Lemma}

\newtheorem{theo}[prop]{Theorem}
\newtheorem{conj}[prop]{Conjecture}
\newtheorem{exam}{Example}[section]

\begin{document}

\baselineskip 16pt

\title{On Spectral Cantor Measures}
\author[I. {\L}aba]{Izabella {\L}aba}
\address{Department of Mathematics  \\University of British Columbia\\ Vancouver, B.C. V6T 1Z2\\Canada}
\email{ilaba@math.ubc.ca}
\author[Y. Wang]{Yang Wang}
\thanks{The second author is supported in part by the National Science Fundation grant
DMS 0070586.}
\address{School of Mathematics  \\ Georgia Institute of Technology \\ Atlanta, Georgia 30332, USA.}
\email{wang@math.gatech.edu}

\date{September 8, 2001}

\begin{abstract}
   A probability measure in $\R^d$ is called a spectral measure if it has an orthonormal basis
consisting of exponentials. In this paper we study spectral Cantor measures. We establish a
large class of such measures, and give a necessary and sufficient condition on the spectrum
of a spectral Cantor measure. These results extend the studies by 
Jorgensen and Pedersen \cite{JoPe98} and Strichartz \cite{Str00}. 
\end{abstract}

\maketitle

\section{Introduction}
\setcounter{equation}{0}

It is known that certain Cantor measures in $\R^d$ have an orthonormal basis consisting
of complex exponentials. This was first observed in Jorgensen and Pedersen \cite{JoPe98}
and studied further in Strichartz \cite{Str00}.  Let $\mu$ be a probability measure in
$\R^d$. We call $\mu$ a {\em  spectral measure} if there exists a $\Lambda\subset\R^d$
such that the set of complex exponentials $\{e(\lambda t):~\lambda\in\Lambda\}$
forms an orthonormal basis for $L^2(\mu)$ (we use $e(t)$ to denote $e^{2\pi it}$
throughout the paper). The set $\Lambda$ is called a {\em spectrum} for $\mu$; we
also say that $(\mu,\Lambda)$ is a {\em spectral pair}.
It should be pointed out that a spectral measure often has more than one spectrum.

In this paper we study spectral Cantor measures in $\R$. Our Cantor measures are 
self-similar measures associated with iterated function systems (IFS).
Consider the iterated functions system (IFS) $\{\phi_j\}_{j=1}^q$ given by
\begin{equation} \label{1.1}
     \phi_j(x) = \rho(x+a_j),
\end{equation}
where $a_j\in\R$ and $|\rho|<1$. It is well known (see e.g. Falconer \cite{Fal-book})
that for any given probability weights $p_1,\dots,p_q >0$ with $\sum_{j=1}^qp_j=1$
there exists a unique probability measure $\mu$ satisfying 
\begin{equation}  \label{1.2}
     \mu = \sum_{j=1}^q p_j \mu\circ \phi_j^{-1}.
\end{equation}
We ask the following question: Under what conditions is $\mu$ a spectral measure?

The familiar middle 3rd Cantor measure given by $\rho=1/3$ and $a_1=0$, $a_1=2$ with $p_1=p_2=1/2$
is not a spectral measure, see Jorgensen and Pedersen \cite{JoPe98}.
The first known example of a spectral measure whose support has non-integer dimension
was given by the same authors in that paper, who showed that
the measure $\mu$ corresponding to $\rho = 1/4$, $a_1=0$, $a_2=1$  and $p_1=p_2=1/2$
is spectral. A spectrum of $\mu$ is
$$
    \Lambda =\Bigl\{ \sum_{k=0}^m d_k 4^k:~m\geq 0,~d_k=0 \mbox{~or~}2\Bigr\}.
$$
Strichartz \cite{Str98} gave an alternative proof of this result, and found 
other examples of spectral Cantor measures. Later Strichartz \cite{Str00}
applied his method to a more general setting to show that a class of
measures having a self-similar type of structure (but not necessarily self-similar
in the more traditional sense) are spectral measures, 
provided that an implicit condition on the zero set of certain trigonometric polynomials 
is satisfied. This condition, however, is not necessary and is somewhat difficult to check.

Spectral measures are a natural generalization of spectral sets. A measurable set $\Omega$
in $\R^d$ with positive and finite measure is called {\em spectral} if $L^2(\Omega)$ has an
orthogonal basis consisting of complex exponentials. Spectral sets have been studied
rather extensively, particularly in recent years. (A partial list of these studies is in
the reference of the paper.) The major unsolved problem concerning spectral sets is the following
conjecture of Fuglede \cite{Fug74}:

\medskip\noindent
{\bf Fuglede's Spectral Set Conjecture:}~ 
Let $\Omega$ be a set in $\R^d$ with positive and finite
Lebesgue measure. Then $\Omega$ is a spectral set if and only if $\Omega$ tiles $\R^d$ by translation.

\medskip
The conjecture remains open in either direction, even in dimension one and for sets
that are unions of unit intervals. As we shall see, the spectral measures--tiling
connection seems to be equally compelling.

In this paper we study self-similar measures satisfying
\begin{equation}  \label{1.3}
     \mu = \sum_{j=1}^q \frac{1}{q} \mu\circ \phi_j^{-1}
\end{equation}
where $\phi_j(x) = \frac{1}{N}(x+d_j)$, $N\in\Z$ and $|N|>1$, and $\D=\{d_j\}\subset\Z$.
We use $\mu_{N,\D}$ to denote the unique probability measure satisfying (\ref{1.3}).
In addition, for each finite subset $\A$ of $\R$ we use $T(N,\A)$ to denote the set
$$
    T(N,\A) :=\Bigl\{\sum_{j=1}^\infty a_j N^{-j}:~a_j\in\A\Bigr\},
$$
which is in fact the attractor of the IFS $\{\phi_a(x):=\frac{1}{N}(x+a):~a\in\A\}$, 
see \cite{Fal-book}. Finally, denote
$$
   \Lambda(N,\A):=\Bigl\{\sum_{j=0}^k a_j N^{j}:~ k\geq 1 \mbox{~and~} a_j\in\A\Bigr\}.
$$

Two finite sets $\A=\{a_j\}$ and $\SS=\{s_j\}$ of cardinality $q$ in $\R$
form a {\em compatible pair}, following the terminology
of \cite{Str00}, if the matrix $M=[\frac{1}{\sqrt{q}}e( a_js_k)]$
is a unitary matrix. In other words $(\delta_{\A}, \SS)$ is a spectral pair, where
$$
    \delta_\A := \sum_{a\in\A}~ \frac{1}{q}\delta(x-a).
$$
For each finite set $\A$ in $\R$ define its {\em symbol} by
$$
     m_\A(\xi):= \frac{1}{|\A|}\sum_{a\in\A} e(-a\xi).
$$
Strichartz \cite{Str00} proves the following theorem:

\begin{theo}[Strichartz] \label{theo-1.1}
    Let $N\in\Z$ with $|N|>1$ and $\D$ be a finite set of integers.
Let $\SS\subset \Z$ such that $0\in\SS$ and
$(\frac{1}{N}\D, \SS)$ is a compatible pair. 
Suppose that $m_{\D/N}(\xi)$ does not vanish on $T(N,\SS)$.
Then the self-similar measure $\mu_{N,\D}$ is a spectral measure with spectrum
$\Lambda(N,\A)$.
\end{theo}
   
Unfortunately the condition that $m_{\D/N}$ does not vanish on $T(N,\SS)$ is not
a necessary condition, and can be very difficult to check, even when both $\D$ and $\SS$ are simple. 
In general we know very little about the zeros $m_{\D/N}$.  Our objective here is
to remove the above condition. We prove that a compatible pair automatically yields
a spectral measure. We also give a necessary and sufficient condition for 
$\Lambda(N,\A)$ to be a spectrum.

\begin{theo} \label{theo-1.2}
Let $N\in\Z$ with $|N|>1$ and let $\D$ be a finite set of integers.
Let $\SS\subset \Z$ such that $0\in\SS$ and
$(\frac{1}{N}\D, \SS)$ is a compatible pair. 
Then the self-similar measure $\mu_{N,\D}$ is a spectral measure. 
If moreover $\gcd(\D-\D)=1$, $0\in \SS$ and $\SS \subseteq [2-|N|, |N|-2]$,
then $\Lambda(N, \SS)$ is a spectrum for $\mu_{N,\D}$.
\end{theo} 
   
In Lemma \ref{lem-2.2} we prove that if $(\frac{1}{N}\D, \SS)$ is a compatible pair
for some $\SS\subset \Z$, then there is also a set $\hat\SS\subset\Z$ 
satisfying the additional conditions of Theorem \ref{theo-1.2} (i.e.,
$0\in \hat\SS$ and $\hat\SS \subseteq [2-|N|, |N|-2]$) such that
$(\frac{1}{N}\D, \hat\SS)$ is a compatible pair.

The following theorem gives a necessary and sufficient condition
for $\Lambda(N,\A)$ to be a spectrum. It also leads to a simple algorithm, see \S3.
 
\begin{theo} \label{theo-1.3}
 Let $N\in\Z$ with $|N|>1$ and $\D\subset \Z$ with $0\in\D$ and 
$\gcd(\D) = 1$.  Let $\SS\subset \Z$ with $0\in\SS$
such that $(\frac{1}{N}\D, \SS)$ is a compatible pair. 
Then $(\mu_{N,\D}, \Lambda(N,\SS))$ is NOT a spectral pair if and only if there exist 
$s_j^*\in\SS$ and nonzero integers $\eta_j$, $0\leq j \leq m-1$, such that 
$\eta_{j+1} = N^{-1}(\eta_j+s_j^*)$ for all $0\leq j \leq m-1$ 
(with $\eta_m :=\eta_0$ and $s_m^*:=s_0^*$).
\end{theo}

The proof of Theorem \ref{theo-1.3} depends on the analysis of the extreme values of the 
eigenfunctions of the Ruelle transfer operator. 
The Ruelle transfer operator was studied in \cite{JoPe98}.

There appears to be a strong link between compatible pairs, tiling of integers
and Fuglede's Conjecture.  All examples suggest that if $\D$ is a finite set of integers
and is part of a compatible pair then $\D$ tiles $\Z$. A finite set $\D\subset \Z$ is
called a {\em complementing set $\tmod{N}$} if there exists a $\E\subset\Z$ such that
$\D \oplus \E$ is a complete residue system $\tmod{N}$. It is known that $\D$ tiles $\Z$
if and only if it is a complementing set $\tmod{N}$ for some $N$. We prove:

\begin{theo}  \label{theo-1.4}
Let $\D\subset\Z$ be a complementing set $\tmod{N}$ with $|N|>1$. Suppose that
$|\D|$ has no more than two distinct prime factors. Then $\mu_{N,\D}$ is a
spectral measure.
\end{theo}

In the next section we shall prove the results just stated. Later in \S3 we give
some examples and state some open problems.

We are indebted to Bob Strichartz for very helpful comments.

\section{Proofs of Theorems}
\setcounter{equation}{0}
 
We first state several lemmas, many of which have been proved in Jorgensen and Pedersen 
\cite{JoPe98} or Strichartz \cite{Str00}. 

\begin{lem}  \label{lem-2.1}
    Let $\A,~\B \subset \R$ be finite sets of the same cardinality. Then the
following are equivalent:
\begin{itemize}
\item[\rm (a)] $(\A,\B)$ is a compatible pair.
\item[\rm (b)] $m_\A(b_1-b_2) = 0$ for any distinct $b_1,\,b_2\in\B$. 
\item[\rm (c)] $\ds \sum_{b\in\B} |m_\A(\xi+b)|^2 \equiv 1$.
\end{itemize}
\end{lem}
\Proof  Note that the condition (b) says precisely that the rows of the matrix
$M=[\frac{1}{\sqrt{|\A|}}e(a_jb_k)]$ are orthonormal. So (a)
and (b) are clearly equivalent.
To see (a) and (c) are equivalent, let
$\delta_\A = \frac{1}{|\A|}\sum_{a\in\A}\delta(x-a)$. Then $\delta_\A$
is a probability measure with $\FT \delta_\A(\xi) = m_\A(\xi)$. Furthermore,
$(\A,\B)$ is a compatible pair if and only if $(\delta_\A,\B)$ is a spectral pair,
see \cite{Str00}. The equivalence of (a) and (c)
follows immediately from Lemma 2.3 of \cite{Str00}. 
\eproof

\begin{lem}  \label{lem-2.2}
    Let $\D,~\SS \subset \Z$ and $N\in\Z$, $|N|>1$ such that $(\frac{1}{N}\D,\SS)$
is a compatible pair. Then
\begin{itemize}
\item[\rm (a)]  $(\frac{1}{N}\D+a, \SS+b)$ is a compatible pair  for any $a, b\in \R$.
\item[\rm (b)]  Suppose that $\widehat S\subset \Z$ such that 
    $\widehat S \equiv S \wmod{N}$. Then $(\frac{1}{N}\D, \widehat S)$ is a
    compatible pair.
\item[\rm (c)]  The elements in both $\D$ and $\SS$ are distinct modulo $N$.
\item[\rm (d)]  Suppose that $|N|>2$. Then there exists an $\widehat\SS$
    with $0\in\widehat\SS$ and $\widehat\SS\subseteq [2-|N|, |N|-2]$
    such that $(\frac{1}{N}\D, \widehat\SS)$ is a compatible pair. 
\item[\rm (e)]  Denote $\D_k = \D+N\D+\cdots+N^{k-1}\D$ and
   $\SS_k = \SS+N\SS+\cdots+N^{k-1}\SS$. Then $(\frac{1}{N^k}\D_k,\SS_k)$
   is a compatible pair.
\end{itemize}
\end{lem}
\Proof  (a) is essentially trivial from Lemma \ref{lem-2.1}. It is also well known
from the fact that any translate of a spectrum is also a spectrum, and any translate
of a spectral measure is also a spectral measure with the same spectra.

For (b), observe that if $s \equiv \hat s \wmod{N}$ then
$m_{\D/N}(\xi + s) = m_{\D/N}(\xi +\hat s)$. Therefore
$$
     \sum_{\hat s\in\widehat \SS}|m_{\D/N}(\xi +\hat s)|^2  
      =\sum_{s\in \SS}|m_{\D/N}(\xi + s)|^2 =1.
$$
This proves (b). 

For (c), assume that $\SS=\{s_j\}$ has $s_1 \equiv s_2 \wmod{N}$. Then we may
replace $s_2$ by $s_1$ in $\SS$ and still have a compatible pair by (b). This
means the matrix $M$ used for defining compatible pairs has two identical
columns, so it cannot be unitary, a contradiction. So elements in $\SS$ are
distinct modulo $N$. Similarly $M$ will have two identical rows if elements
in $\D$ are not distinct modulo $N$, again a contradiction.

To prove (d), we first translate $\SS$ so that $0\in\SS_1:=\SS+a$ for some $a\in\Z$.
$(\frac{1}{N}\D,\SS_1)$ is still a compatible pair. Now $[2-|N|, |N|-2]$
contains a complete set of residues $\tmod{N}$ because it contains
at least $|N|$ consecutive integers. Choose $\widehat\SS \subseteq
[2-|N|, |N|-2]$ so that $0\in\widehat\SS$ and 
$\widehat\SS\equiv \SS_1 \wmod{N}$. Then $(\frac{1}{N}\D,\widehat\SS)$ is
a compatible pair.

Finally (e) is a special case of Lemma 2.5 in \cite{Str00}.
\eproof

\begin{lem}  \label{lem-2.3}
   Under the assumptions of Theorem \ref{theo-1.2}, let 
$Q(\xi):=\sum_{\lambda\in\Lambda(N,\SS)}|\FT\mu(\xi +\lambda)|^2$. 
Then
\begin{itemize}
\item[\rm (a)] The set of exponentials $\{e(\lambda \xi):~\lambda\in\Lambda(N,\SS)\}$ 
           is orthonormal in $L^2(\mu)$.
\item[\rm (b)] $Q(\xi) \leq 1$ for all $\xi\in\R$ and 
           $\{e(\lambda \xi):~\lambda\in\Lambda(N,\SS)\}$
          is an orthonormal basis for $L^2(\mu)$ if and only if $Q(\xi)\equiv 1$.
\item[\rm (c)] The function $Q(\xi)$ is the restriction of an entire function of
          exponential type to the real line. Furthermore it satisfies
         \begin{equation}  \label{2.1}
           Q(\xi) = \sum_{s\in\SS} |m_\D(N^{-1}(\xi+s))|^2 Q(N^{-1}(\xi+s)).
         \end{equation}
\end{itemize}
\end{lem}
\Proof  See Jorgensen and Pedersen \cite{JoPe98}. The right hand side of (\ref{2.1}) is known as the {\em Ruelle transfer operator} (operated on $Q$).
\eproof

\bigskip

\noindent
{\bf Proof of Theorem \ref{theo-1.3}.}~~Denote $\Lambda :=\Lambda(N,\SS)$ and $\mu:=\mu_{N,\D}$. 

\vspace{2mm}
\noindent
($\Leftarrow$)~~We prove that $(\mu, \Lambda)$ is not a spectral pair 
by proving that $Q(\xi) \not\equiv 1$, where $Q(\xi)$ is defined in Lemma~\ref{lem-2.3}. 
In fact, we prove that $Q(\eta_0) = 0$.

Observe that $m_\D(\eta_j) =1$ because $\eta_j\in\Z$. Now 
$\sum_{s\in\SS} |m_\D(N^{-1}(\eta_j+s))|^2  =1$ by (c) of Lemma \ref{lem-2.1}. 
Since $m_\D(N^{-1}(\eta_j+s_j^*))=m_\D(\eta_{j+1}) =1$, 
it follows that $m_\D(N^{-1}(\eta_j+s))=0$ for all $s\neq s_j^*$ in $\SS$.

Take any $\lambda\in\Lambda$ and write $\lambda= \sum_{k=0}^\infty s_k N^k$ where $s_k\in\SS$ and of course only finitely many $s_k \neq 0$. We have
$$
      \FT\mu(\eta_0+\lambda) = \prod_{j=1}^\infty m_\D(N^{-j}(\eta_0+\lambda)).
$$
Note that 
$$
     m_\D(N^{-1}(\eta_0+\lambda))= 
                  m_\D\Bigl(N^{-1}(\eta_0+s_0)+\sum_{k=0}^\infty s_{k+1}N^k\Bigr)
                 = m_\D(N^{-1}(\eta_0+s_0)).
$$
Hence $m_\D(N^{-1}(\eta_0+\lambda))= 1$ for $s_0=s_0^*$ and $m_\D(N^{-1}(\eta_0+\lambda))= 0$ otherwise. Suppose $s_0=s_0^*$. Then the same argument together with
the fact $\eta_1 = N^{-1}(\eta_0+s_0^*)$ yield $m_\D(N^{-2}(\eta_0+\lambda))= 1$ 
for $s_1=s_1^*$ and 
$m_\D(N^{-2}(\eta_0+\lambda))= 0$ otherwise. By induction we easily obtain
$m_\D(N^{-j}(\eta_0+\lambda)) \neq 0$ if and only if $s_j = s^*_{j\wmod{m}}$. Therefore
$\FT\mu(\eta_0+\lambda) \neq 0$ only if $s_j = s^*_{j\wmod{m}}$ for all $j \geq 0$. But this is
impossible since $s_j = 0$ for all sufficiently large $j$. 
Thus $\FT\mu(\eta_0+\lambda)=0$ and $Q(\eta_0)=0$. 

\vspace{2mm}

\noindent
($\Rightarrow$)~~Assume that  $(\mu, \Lambda)$ 
is NOT a spectral pair. Then $Q(\xi)\not \equiv 1$. Note that this can not happen if
$|\D|=|\SS|=1$, in which case
$\mu = \delta_0$ and $\Lambda = \{0\}$. So we may assume that
$q=|\D|=|\SS|>1$. It is well known that in this case $T:=T(N,\SS)$ is a
compact set with infinite cardinality. 

Since $Q(\xi) \not\equiv 1$, $Q(\xi) \not \equiv 1$ for $\xi\in T$, because
$Q$ is extendable to an entire function on the complex plane and $T$ is an infinite compact set. 
Denote $X^-:=\{\xi\in T:~ Q(\xi) = \min_{\eta\in T}Q(\eta)\}$. It follows from $Q(0) =1$ that
$0\not\in X^-$.We apply the Ruelle transfer operator to derive a contradiction.

   For any $s\in\R$ denote $\phi_s(\xi) = N^{-1}(\xi+s)$. Then
$T = \bigcup_{s\in\SS}\phi_s(T)$. Hence $\phi_s(\xi) \in T$ for all $\xi\in T$.
Now choose any $\xi_0\in X^-$ and set $Y_0 = \{\xi_0\}$. Define recursively
$$
      Y_{n+1} = \bigl\{\phi_s(\xi):~s\in\SS,\,\xi\in Y_n,\,
         \phi_s(\xi)\in X^-\bigr\}
      \mhsp \mbox{(counting multiplicity)}.
$$

\noindent
{\it Claim 1: We have $|Y_{n+1}| \geq |Y_n|$ (counting multiplicity).}
\vspace{1ex}

\noindent
{\it Proof of Claim 1.}  Let $\xi^*\in X^-$. By (c) of Lemma \ref{lem-2.3}
$$
      \min_{\eta\in T} Q(\eta) = Q(\xi^*) = 
          \sum_{s\in\SS} |m_\D(\phi_s(\xi^*))|^2 Q(\phi_s(\xi^*)).
$$
But $\sum_{s\in\SS} |m_\D(\phi_s(\xi^*))|^2=1$ by Lemma \ref{lem-2.3}
and $Q(\phi_s(\xi^*))\geq \min_{\eta\in T}Q(\eta)$. Thus
$Q(\phi_s(\xi^*))= \min_{\eta\in T}Q(\eta)$ whenever  $m_\D(\phi_s(\xi^*))\neq 0$.
In other words,
\begin{equation}  \label{2.2}
     \phi_s(\xi^*)\in X^{-} \shsp\mbox{whenever} \shsp m_\D(\phi_s(\xi^*))\neq 0.
\end{equation}
Hence for each $\xi\in Y_n$ there exists at least one $s\in\SS$ such that
$\phi_s(\xi)\in Y_{n+1}$, proving Claim 1.
\esubproof

\vspace{1ex}

\noindent
{\it Claim 2: The elements of $Y_n$ in fact all have multiplicity one.}
\vspace{1ex}

\noindent
{\it Proof of Claim 2.}  It is easy to see that elements in $Y_n$ have the form
$\phi_{s_n}\circ\cdots\circ\phi_{s_1}(\xi_0)$. If some element in $Y_n$ has multiplicity
more than one, then there are two distinct sequences 
$(s_1,\dots, s_n)$ and $(t_1, \dots, t_n)$ in $\SS$ such that
$$
    \phi_{s_n}\circ\cdots\circ\phi_{s_1}(\xi_0) = 
     \phi_{t_n}\circ\cdots\circ\phi_{t_1}(\xi_0).
$$
Expanding the two expressions yields
$$
    \frac{1}{N^n}(\xi_0+s_1+Ns_2 +\cdots + N^{n-1}s_n) =
    \frac{1}{N^n}(\xi_0+t_1+Nt_2 +\cdots + N^{n-1}t_n).
$$
But this is clearly not possible, since all elements of $\SS$ are in different
residue classes $\tmod{N}$.
\esubproof

Now $X^-$ is finite. It follows that all $Y_n$ have the same cardinality for
sufficiently large $n$, say $n \geq n_0$. Therefore for $n\geq n_0$ any
$\xi^*\in Y_n$ has a unique offspring $\phi_s(\xi^*)\in Y_{n+1}$. Furthermore,
for any $t\in\SS$ and $t\neq s$ we must have $m_{\D}(\phi_t(\xi^*)) = 0$, 
because otherwise by (\ref{2.2}) we would have 
$\phi_t(\xi^*)\in Y_{n+1}$ as another offspring, a contradiction. 
Thus starting with a $\xi_{n_0}\in Y_{n_0}$ we obtain a
sequence $\{\xi_n\}_{n\geq n_0}$ in which $\xi_n\in Y_n$ 
is the unique offspring of $\xi_{n-1}\in Y_{n-1}$, $\xi_n = \phi_{s_{n-1}}(\xi_{n-1})$
for some $s_{n-1}\in \SS$.
It follows from the finiteness of $X^-$ that there exist $k>n_0$ and $m>0$ such that
$$
     \xi_k = \xi_{k+m} = \frac{1}{N^m}(\xi_k+ s_k+s_{k+1}N+\cdots+s_{k+m}N^{m-1}).
$$
In particular $\xi_k\in\Q$. Set $\eta_j = \xi_{j+k}$ and $s_j^* = s_{j+k}$ 
for $0 \leq j \leq m-1$
($\eta_m :=\eta_0$ and $s_m^*:=s_0^*$. Only the case $j=m-1$ needs to be check).
. Then $\eta_j\in\Q$ for all $j$. Furthermore 
$\eta_{j+1} = \phi_{s_j^*}(\eta_j) = N^{-1}(\eta_j+s_j^*)$ for all $0\leq j \leq m-1$

Note that $m_{\D}(\phi_s(\eta_j)) = 0$ 
for all $0\leq j \leq m-1$ and $s \neq s_{j}^*$, as shown above. This yields 
$$
       |m_\D(\eta_{j})|^2 = \sum_{s\in\SS} |m_\D(\phi_{s_{j-1}^*}(\eta_{j-1}))|^2  =1, \mhsp
        1\leq j \leq m.
$$
But $m_\D(\eta_j) = \frac{1}{|\D|}\sum_{d\in\D} e(d \eta_j)$ and $0\in\D$. The only
way $|m_\D(\eta_j)|=1$ can hold is $e(d \cdot\eta_j) = 1$ for all $d\in\D$.
Hence $\eta_j\cdot d \in \Z$ for all $d\in\D$. Since all $\eta_j \neq 0$
and $\gcd(\D)=1$, it follows that all $\eta_j \in\Z$. The theorem is proved.
\eproof

\vspace{1ex}

\noindent
{\bf Proof of Theorem~\ref{theo-1.2}.}~~Suppose that $\mu_{N,\D}$ is a spectral
measure with spectrum $\Lambda$ then $\mu_{N, a\D}$ is a spectral measure
with spectrum $a^{-1}\Lambda$. Therefore to prove that $\mu_{N,\D}$ is spectral
we may without loss of generality assume that $0\in\D$ and $\gcd(\D) = 1$.

If $|N|=2$ then $\D=\{0,1\}$ or $\D=\{0,-1\}$. The corresponding self-similar measure
$\mu_{N,\D}$ is the Lebesgue measure supported on a unit interval. Clearly
in this case $\mu_{N,\D}$ is a spectral measure with spectrum $\Z$.

For $|N|\geq 3$, by Lemma \ref{lem-2.2} we may replace $\SS$ by
$\widehat\SS \subseteq [2-|N|,|N|-2]$ such that $0\in\widehat \SS$ and
$(\frac{1}{N}\D,\widehat\SS)$ is a compatible pair.
Now for any such $\widehat\SS$, we have
$T(N,\widehat\SS)\subseteq [-\frac{N-2}{N-1}, \frac{N-2}{N-1}]$ for $N>0$
or $T(N,\widehat\SS)\subseteq [-\frac{N^2+N-2}{N^2-1},\frac{N^2+N-2}{N^2-1}]$
for $N<0$. In either case $T(N,\widehat\SS)$ contains no integer other than 0.
Now suppose that $(\mu_{N,\D}, \Lambda(N,\widehat\SS))$ is not a spectral pair. Then
there exist $s_j^*$ in $\SS$ and nonzero integers $\eta_j$ satisfying the condition
of Theorem \ref{theo-1.3}. Starting with $T_0 = \{\eta_j\}$ we see that
$T_0 \subseteq \bigcup_{s\in\SS}\phi_s(T_0)$. This yields $T_0 \subseteq T(N,\SS)$.
But this is a contradiction since $T(N,\widehat\SS)$ contains no integer other than 0.
Therefore $(\mu_{N,\D}, \Lambda(N,\widehat\SS))$ is a spectral pair by Theorem
\ref{theo-1.3}.
\eproof

\bigskip

\noindent{\bf Proof of Theorem \ref{theo-1.4}.}~~It is known that if $\D$ is a 
complementing set $\tmod {N}$ then there
exists an $L|N$ whose prime factors are precisely those of $|\D|$ such that
$\D$ is a complementing set $\tmod {L}$. We prove that there exists an
$\SS\in\Z$ such that $(\frac{1}{N}\D,\SS)$ is a compatible pair, using a
theorem of Coven and Meyerowitz \cite{CoMe99}.  The argument below
is essentially a repetition of the proof of Theorem 1.5(i) in \cite{Lab02}
specialized to the case of two prime factors.

Let $\Phi_n(z)$ denote the $n$-th cyclotomic polynomial.  Let also
$D(z)=\sum_{d\in\D} z^d$ so that $m_\D(\xi)=D(e^{2\pi i\xi})$.
Assume that $|\D| = p^\alpha q^\beta$ and
$L = p^{\alpha'} q^{\beta'}$, where $p,\,q$ are distinct primes. (If
$|\D|$ is a prime power, the proof below works and is simpler.) Let
$$
    \P :=\Bigl\{p^k:~\Phi_{p^k}(z)~|~D(z),~k\leq \alpha'\Bigr\}, \shsp
    \QQ :=\Bigl\{q^k:~\Phi_{q^k}(z)~|~D(z),~k\leq \beta'\Bigr\}.
$$
Coven and Meyerowitz \cite{CoMe99} prove that $|\P| = \alpha$, $|\QQ| = \beta$,
and that $\Phi_{p^kq^l}(z)~|~m_\D(z)$ for all $k,l$ such that $p^k\in\P,
q^l\in\QQ$.

We construct the set $\SS$. Write $\P=\{p^{k_j}\}$ and $\QQ=\{q^{l_j}\}$, where
$k_1<k_2< \cdots <k_\alpha$ and $l_1<l_2<\cdots l_\beta$. Define
$$
   \E =\Bigl\{\sum_{j=1}^\alpha a_j p^{-k_j} +\sum_{j=1}^\beta b_j q^{-l_j}:~
                0 \leq a_j <p,~0\leq b_j <q\Bigr\},
$$
and let $\SS = N\E$. Clearly $\SS\subset\Z$.
To prove that $(\frac{1}{N}\D,\SS)$ is a compatible pair it suffices to prove that
$m_{\D/N}(s_1-s_2) = 0$ for any distinct $s_1, s_2\in\SS$. Equivalently, we only need
to show that $m_\D(\lambda_1-\lambda_2) = 0$ for any distinct $\lambda_1, \lambda_2
\in\E$. Note that 
$$
      \lambda_1-\lambda_2 =\sum_{j=1}^\alpha c_j p^{-k_j} +\sum_{j=1}^\beta d_j q^{-l_j},
      \shsp  -p <c_j <p,~-q<d_j< q.
$$
If all $c_j =0$ then not all $d_j =0$. So $\lambda_1-\lambda_2 = r/q^{l_{j^*}}$
for some $r$ with $\gcd(r, q) = 1$, where $j^*$ is the largest $j$ such that $d_j \neq 0$.
Therefore $\Phi_{q^{l_{j^*}}}(e^{2\pi i(\lambda_1-\lambda_2)}) = 0$ and hence
$m_\D(\lambda_1-\lambda_2) = 0$. Similarly $m_\D(\lambda_1-\lambda_2) = 0$ if all
$d_j = 0$. Finally, assume none of the above is true. Let $j_1$ be the largest $j$
such that $c_j \neq 0$ and $j_2$ be the largest $j$ such that $d_j \neq 0$.
Then $\lambda_1-\lambda_2 = r_1/p^{k_{j_1}} + r_2/q^{l_{j_2}}$ with
$\gcd(r_1, p) = \gcd(r_2, q) = 1$. This yields
$$
      \lambda_1-\lambda_2 = \frac{r_1 q^{l_{j_2}} +r_2 q^{k_{j_1}}}{p^{k_{j_1}}q^{l_{j_2}}}.
$$
The numerator is clearly coprime with the denominator. 
Hence $\Phi_{p^{k_{j_1}}q^{l_{j_2}}}(e^{2\pi i(\lambda_1-\lambda_2)})=0$ and therefore
$m_\D((\lambda_1-\lambda_2)=0$.

So $(\D/N, \SS)$ is a compatible pair, and hence $\mu_{N,\D}$ is a spectral measure.
\eproof

\section{Examples and Open Questions}
\setcounter{equation}{0}

Theorem \ref{theo-1.3} leads to an algorithm for determining whether $\Lambda(N,\SS)$
is a spectrum for $\mu_{N,\D}$. To find whether the integer sequence
$\{\eta_j\}$ exists we only need to check a finite number of integers.
This is because we have shown in the proof of Theorem \ref{theo-1.3} that if
the sequence $\{\eta_j\}$ exists, it must be contained
in $T(N,\SS)$. However, $T(N,\SS)$ is compact. In fact 
$T(N,\SS)\subseteq [\frac{a}{N-1}, \frac{b}{N-1}]$ for $N>0$ and
$T(N,\SS)\subseteq [\frac{a+Nb}{N^2-1}, \frac{b+Na}{N^2-1}]$ for $N<0$,
where $a, b$ are the smallest and the largest elements in $\SS$, respectively.
\vspace{1ex}

\begin{exam} {\rm Our first example addresses the condition in Theorem ~\ref{theo-1.1}
by Strichartz. Let $N=5$ and $\D=\{0,\pm 2, \pm 11\}$. Since $\D$ is a residue system 
$\tmod{N}$, the set $T(N,\D)$ is a fundamental domain of $\Z$ and $\mu_{N,\D}$ is 
simply the restriction of the Lebesgue measure to $T(N,\D)$, see Gr\"ochenig and Haas \cite{GrHa94}.
Let $\SS=\{0,\pm 1,\pm 2\}$. Then $\Lambda(N,\SS) = \Z$ is a spectrum for $\mu_{N,\D}$.
However, $m_\D(\xi)$ has a zero in $[0, \frac{1}{4}]$, which is contained in
$T(N,\SS) = [-\frac{1}{2}, \frac{1}{2}]$.

The condition can also be hard to check. Consider the same $\D$ as above, let $N=10$ and
$\SS=\{0,2,4,6,8\}$. Then $(\D/N,\SS)$ is a compatible pair. By Theorem \ref{theo-1.2}
$\Lambda(N,\SS)$ is a spectrum for $\mu_{N,\D}$. Nevertheless it is difficult to check
whether $T(N,\SS)$ contains a zero of $m_\D(\xi)$.}
\eproof
\end{exam}

\vspace{1ex}

\begin{exam}  {\rm Our next example illustrates that a spectral measure can have many spectra.
Take $N=6$ and $\D=\{0,1,2\}$. Then $(\D/N, \SS)$ are compatible pairs for
both $\SS=\{0,2,4\}$ or $\SS=\{0,-2,2\}$. By Theorem \ref{theo-1.2} $\Lambda(N,\SS)$
are spectra of $\mu_{N,\D}$ for both $\SS$.

    A far more striking example is to take $\SS'=\{0,4,8\}$. Then $(\D/N,\SS')$ is a
compatible pair because $\SS' \equiv \SS=\{0,2,4\} \wmod{N}$. One can check using the
algorithm described earlier that $\Lambda(N,\SS')$ is indeed a spectrum of $\mu_{N,\D}$,
as is $\Lambda(N,\SS)$. But $\Lambda(N,\SS') = 2\Lambda(N,\SS)$! This is rather striking
because $\Lambda(N,\SS')$ is intuitively ``twice as sparse as'' $\Lambda(N,\SS)$.}
\eproof
\end{exam}

The study in this paper also leaves several questions unanswered. For example, do Theorem
\ref{theo-1.2} and Theorem \ref{theo-1.3}, or something similar, hold in higher dimensions? 
The difficulty is that an analytic function of two
or more variables may attain its infimum at infinitely many points, even on a compact set.
The technique of the Ruelle transfer operator employed in this paper has its
origin in the study of wavelets and self-affine tiles, see e.g. \cite{GrHa94}, in 
dimension 1. It was extended to higher dimensions in Lagarias and Wang \cite{LaWa00}.
Could the techniques there be applied to higher dimensions to yield results on spectral measures?

Note that we have studied spectral Cantor measures in which the probability
weights are equally distributed. Is this a general rule? We conclude this paper with the following conjecture:

\begin{conj} Let $\mu$ be the self-similar measure associated with the IFS 
$\phi_j(x)= \rho(x+a_j)$, $1\leq j \leq q$, with probability weights $p_1, \dots,p_q>0$,
where $|\rho|<1$. Suppose that $\mu$ is a spectral measure. Then
\begin{itemize}
\item[\rm (a)]  $\rho = \frac{1}{N}$ for some $N\in\Z$.
\item[\rm (b)]  $p_1=\cdots=p_q = \frac{1}{q}$.
\item[\rm (c)]  Suppose that $0\in {\mathcal A} =\{a_j\}$. Then $\A = \alpha \D$
for some $\alpha\in\R$ and $\D\subset\Z$. Furthermore, $\D$ must be a 
complementing set $\tmod{N}$.
\end{itemize}
\end{conj}


\end{document}